\documentclass[12pt]{article}
\usepackage{amsmath,amsthm,amssymb}
\usepackage{graphicx}
\begin{document}
\title{Chebyshev and garment cutting.\\ Debunking some myths}

\author{Victor Tapia\\
Via dei Giardini 6, Trieste, 34146, Italy\\
tapiens@gmail.com}

\maketitle

\abstract{In {\tt 1878}, Pafnuty Chebyshev presented to the {\it Association fran\c caise pour l'avan\-cement des sciences} {\it [French Association for the Advancement of the Sciences]} an article \cite{Chebyshev1878} dealing with garment cutting. According to Chebyshev himself, his interest was sparked by a lecture given by \'Edouard Lucas that he had attended in {\tt 1876} \cite{Lucas1876}.

There is a second story on the origin of Chebyshev's interest in garment cutting according to which in the 1850s, being short of money, Chebyshev got himself a job as a consultant to a clothing factory. At the time of the Crimean War (1853-1856), there was a great demand for uniforms. Chebyshev was allegedly asked to optimize the use of fabric, and it was there that his interest in garment cutting was born.

This second story appears to have its origin in a post by Clive J. Grant to MacTutor in 1996 \cite{Grant1996}. However, this contribution contains no references, and no other source of information that I have found offers any first-hand documentation to support this story. Our conclusion is that this second story is a fabrication, invented out of whole cloth.}

\section*{The Facts}

Pafnuty Lvovich Chebyshev was born in Russia in 1821 near Moscow. In 1837, he entered Moscow University, and in 1847, upon graduation, he moved to Saint Petersburg, where he became an adjunct professor of mathematics at Saint Petersburg University. In 1850, he was promoted to the rank of associate professor, and to full professor in 1860. Chebyshev was in contact with Western European mathematicians, mainly in France, where he attended several congresses.

Chebyshev's mathematical output was enormous; see \cite{Chebyshev1946b}. Here we mention only those publications relevant to our concern. In 1876, Chebyshev attended a lecture by \'Edouard Lucas \cite{Lucas1876} on garment cutting at the {\it 5me session de l'Association fran\c caise pour l'avancement des sciences} in Clermont-Ferrand. The proceedings of the session contain only the title of that lecture.

In 1878, the seventh session of the Association fran\c caise pour l'avancement des sciences took place in Paris. Chebyshev presented four communications; one of them \cite{Chebyshev1878} is entitled {\it Sur la coupe des v\^etements} {\it [On the cutting of garments]}. A brief account, one page long, appeared in the proceedings of the session. As Chebyshev himself recognized in the report of that lecture, the motivation for his interest was the above--mentioned lecture of 1876 by Lucas \cite{Lucas1876}. Chebyshev did not publish the full text elsewhere, as he sometimes did in similar cases.

Chebyshev's article is followed, in the same issue of the proceedings, by a half--page abstract of a contribution by Lucas \cite{Lucas1878}, apparently not from his hand (Figure 1). In 1880, Lucas published, in Italian, a longer article, which was translated into French and published in 1911 \cite{Lucas1911}.

Before his death, Chebyshev indicated in his will that only those works that he himself had marked ``printable'' should be published. There was no such mark on the article ``On the cutting of garments.”

The manuscript of Chebyshev's article was found only after his death in 1894. Since that manuscript did not bear the word ``printable,'' Markoff and Sonin \cite{Markoff1899} did not include it in their edition of the {\it \OE uvres de P. L. Tchebychef} (1899, 1907), but they provided a brief comment similar to what had appeared in the proceedings of the {\it Association fran\c caise pour l'avancement des sciences} in 1878. A different choice was made in the {\it Complete Collected Works} \cite{Chebyshev1946b}, published in 1946, in which we find a full Russian translation of the manuscript. In some sense, Chebyshev's will was respected, since the original version in French was not published. That same year, a translation into Russian \cite{Chebyshev1946a} was independently published by Popov. A translation into English was published by Chobot and Collomb in 1970 \cite{Chobot1970}.

This account is based on several credible sources, among them Vassilief (1898) \cite{Vassilief1898}, Poss\'e in Markoff et Sonin's {\it \OE uvres} (1899, 1907) \cite{Markoff1899}, Chobot and Collomb (1970) \cite{Chobot1970}, Youschkevitch in {\it Dictionary of Scientific Biography} (1970-1980) \cite{Youschkevitch1970}, and Butzer and Jongmans (1999) \cite{Butzer1999}.

A good summary of the story can be found in {\it \OE uvres} \cite{Markoff1899}:

\begin{itemize}

\item[]{\it Apr\`es avoir indiqu\'e que l'id\'ee de cette \'etude lui est venue lors de la communication faite, il y a deux ans, au Congr\`es de Clermont-Ferrand, par M. \'Edouard Lucas, sur la g\'eom\'etrie du tissage des \'etoffes à fils rectilignes, M. Tchebichef pose les principes g\'en\'eraux pour d\'eterminer les courbes suivant lesquelles on doit couper les diff\'erents morceaux d'une \'etoffe, pour en faire une gaine bien ajust\'ee, servant à envelopper un corps de forme quelconque.}

\item[]{\it After indicating that the idea for this study came to him during the communication made two years before at the Clermont-Ferrand Congress, by Mr. \'Edouard Lucas, on the geometry of weaving fabrics with straight threads, Mr. Chebyshev sets out the general principles for determining the curves along which we must cut the different pieces of fabric to make a well-fitting sheath, used to wrap a body of any shape.}

\end{itemize}

All these first-hand sources indicate that the origin of Chebyshev's interest in garment cutting was the lecture by Lucas in 1876, as he himself declared.

\bigskip


\begin{center}

\includegraphics[scale=1.0, height=12.4truecm, width=14.8truecm]{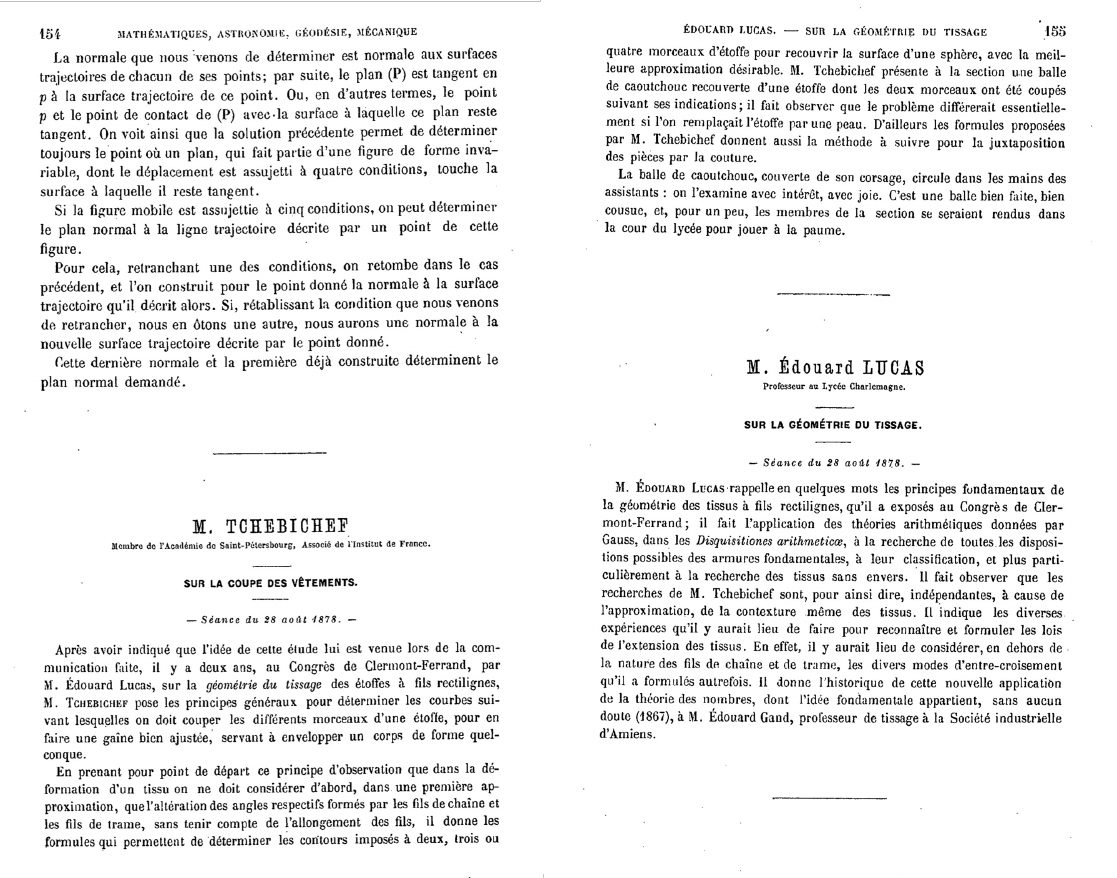}

\end{center}

\small

\noindent {\bf Figure 1.} Two pages of {\it Compte Rendu de la 7 e session de l'Association fran\c caise pour l'avancement des sciences} containing the reports of the lectures by Chebyshev and Lucas. (Reproduced by kind permission of the {\it Association fran\c caise pour l'avancement des sciences} (AFAS).)

\normalsize


\bigskip

During the eighteenth and nineteenth centuries, cartography was a fashionable subject among mathematicians. Its mathematical study was initiated by Euler in 1777. Later, in 1779, Lagrange \cite{Lagrange1779} published his {\it Sur la construction des cartes g\'eographiques} {\it [On the construction of geographic maps]}.

In 1856, Chebyshev published two articles \cite{Chebyshev1856a,Chebyshev1856b} on cartography with the same title, {\it Sur la construction des cartes g\'eographiques}. He was heavily influenced by Lagrange and used the latter's work in his own articles on cartography. Even more, the title of Chebyshev's two papers on cartography is the same as that of the articles by Lagrange.

A second important influence on Chebyshev was the development of analytic geometry. As a consequence, the problem of the cartographic projection of geographical regions with minimal deformation took a new direction. One important contribution is due to Gauss \cite{Gauss1822}, in 1822, with the suggestive title {\it General solution of the problem to represent the parts of a given surface on another given surface so that the smallest parts of the representation shall be similar to the corresponding parts to the surface represented}. For a detailed historical account, see \cite{Papadopoulos2022a,Papadopoulos2022b}.

Cartography is relevant to our story because as mathematical problems, cartography and garment cutting are related. Therefore, a knowledge of cartography is almost a must for the problem of garment cutting. Chebyshev developed that knowledge in his cartography publications of 1856.

In 1908, Gaston Darboux \cite{Darboux1909}, in his contribution to the International Congress of Mathematicians in Rome, mentions questions related to geographical maps and also mentions the work of Chebyshev on garment cutting. He includes both subjects in his paper {\it Les origines, les m\'ethodes et les probl\`emes de la g\'eom\'etrie infinit\'esimale} {\it [The origins, methods and problems of infinitesimal geometry]}. The similarity is suggested by the title of the 1822 article by Gauss on conformal mappings. The mathematical problem was brilliantly summarized by Athanase Papadopoulos in 2016 \cite{Papadopoulos2016}:

\begin{itemize}

\item[]{\it The problem of drawing geographical maps and the one of fitting of garments are inverse of each other. Indeed, on the one hand, one searches for mappings from a piece of the sphere into a Euclidean piece of paper with minimal distortion, and on the other hand, one constructs a map from a Euclidean piece of fabric onto a curved surface (part of a human body), such that the fabric fits the curved surface with minimal distortion.}

\end{itemize}

Chebyshev's interest in garment cutting was posterior to his work on cartography. Here we can only speculate that his interest in garment cutting was facilitated by his knowledge of cartography.

\section*{Fiction}

In 1996, Clive J. Grant \cite{Grant1996} posted a contribution on MacTutor entitled {\it Chebyshev Nets}, which we reproduce here in full:

\begin{itemize}

\item[]{\it What nobody writes about is something you may wish to consider as a useful topic for inclusion in a history of mathematics and their application to engineering problems. (Most engineers don't have the foggiest idea of who Chebyshev was, even from his mathematical works.)}

\item[]{\it Not unlike some of his modern counterparts, Chebyshev was paid a pittance by his university. To supplement his income, he took on private clients. Among the latter was the owner of a huge textile business. With the onset of the Crimean War, there was a huge demand for army uniforms. Chebyshev was assigned the task of developing a means of cutting fabrics more economically. All concerned seemed to believe that the difficulty of forming the shoulder seam should be the first task that Chebyshev took on. As might be expected, he saw the fabric, not as a limp piece of material, but as a kind of net. He became so interested in this that he abandoned his client and went on to develop a very important theory of cable nets. He wrote a book on the subject, but it has never been translated. In addition to the arcane mathematical notations of mid-nineteenth century Russia, the topic, per se, is of little interest except to a very narrow group of engineers. I think there may be fifteen of us all told. The class of cablenet structures Chebyshev described are now referred to as {\sc Chebyshev nets}.}

\end{itemize}

This contribution contains no references, and this version of the story has been repeated with different wording in several subsequent works. In order to help trace the transmission of this version of the story, we have added the citations on which each article we cite is based.

\begin{itemize}

\item[]Fuchs and Tabachnikov (2007) \cite{Fuchs2007}:

\begin{itemize}

\item[]{\it Pafnuty Chebyshev, a prominent Russian mathematician of the 19th century, was motivated by an applied problem: how to cut fabric more economically (he was working for a private client, an owner of a textile business). This was an acute problem: with the onset of the Crimean War, there was a huge demand for army uniforms.}

\item[$\bullet$]This account is similar to Grant's, but no reference is given, either to Grant or to any other source.

\end{itemize}

\item[]Ghys (2011) \cite{Ghys2011}:

\begin{itemize}

\item[]{\it Selon Grant [22], pour compl\'eter son maigre salaire universitaire, Tche\-bychev avait obtenu un contrat pour optimiser la d\'ecoupe d'uniformes militaires \'evitant les pertes de tissu. Le probl\`eme \'etait semble-t-il particuli\`erement d\'elicat pour les \'epaules $\cdots$. Cette question l'int\'eressa tant qu'il en oublia son client pour d\'evelopper une th\'eorie des r\'eseaux qu'on appelle aujourd'hui ``de Tchebychev.'' Il aurait \'ecrit un livre sur la question, que je n'ai malheureusement pas pu localiser.}

\item[]

\item[$\bullet$]The reference is to Grant's posting \cite{Grant1996}.

\end{itemize}

\item[]The review by Athanase Papadopoulos (2013) \cite{Papadopoulos2013}:

\begin{itemize}

\item[]{\it The first few lines of this paper, as well as a relatively long appendix, contain some history of the subject. One learns that Chebyshev, in order to supplement his income as a university professor, signed a contract with a textile businessman, for a second job where the goal was to find methods of cutting military uniforms (for which there was a huge demand at that time) that would minimize fabric loss. This was Chebyshev's motivation for studying the cutting of garments. It is not clear to what extent Chebyshev helped the profession of clothing design, because he immediately got more interested in the mathematical aspect of the theory. This led Chebyshev to write a paper which bears the title of the paper under review, and a book which seems to be very difficult to locate. Chebyshev's paper, written in French, is reproduced in the appendix of the paper under review. In particular, Chebyshev developed the theory of what are now called {\sc Chebyshev nets}.}

\item[$\bullet$]The review cites Ghys (2011) \cite{Ghys2011}.

\end{itemize}

\item[]Papadopoulos (2016) \cite{Papadopoulos2016}:

\begin{itemize}

\item[]The Grant version is not quoted and the 2011 article by Ghys \cite{Ghys2011} is cited but is not essential to the article.

\end{itemize}

\item[]Khesin and Tabachnikov (2021) \cite{Khesin2021}:

\begin{itemize}

\item[]{\it Pafnuty Chebyshev $\cdots$ came up with this notion, motivated by an applied problem: how to cut fabric in a more economical way (there was a high demand for army uniforms due to the Crimean War). Chebyshev presented his results in a talk entitled {\sc Sur la coupe des v\^etements} in {\tt 1878}, at a session of the {\sc Association fran\c caise pour l'avance\-ment des sciences in Paris}.}

\item[$\bullet$]Cited are Fuchs and Tabachnikov (2007) \cite{Fuchs2007}, Ghys (2011) \cite{Ghys2011}.

\end{itemize}

\item[]Papadopoulos (2021) \cite{Papadopoulos2021}:

\begin{itemize}

\item[]{\it According to this report, the idea for this work occurred to him during a communication made two years earlier by Lucas at the meeting of the same association in Clermont-Ferrand.}

\item[]{\it The following text is contained in the MacTutor website of History of Mathematics, at the article Chebyshev net(s) and it is attributed to Clive J. Grant. The story might not be authentic, but we quote it because it gives a lively idea of the general atmosphere in which Chebyshev lived.}

\item[$\bullet$]Grant (1996) \cite{Grant1996} is quoted in the text but not cited in the references.

\end{itemize}

\end{itemize}

\section*{Grant's Account}

The account given by Grant departs in several respects from the historical record. Some important points of departure are Chebyshev's lack of money, the connection between uniforms and the army, and the existence of an untranslated book. Let us analyze the statements made by Grant.

\bigskip

{\bf Money.} From the biography contained in his {\it \OE uvres}, volume II (1907), we learn that when Chebyshev settled in Saint Petersburg in 1847, his financial situation was critical. But this shortage of money is unrelated to our story, which is supposed to have taken place during the Crimean War, that is, between 1853 and 1856. Chebyshev's parents had once been wealthy, but they had lost their lands a few years previously, in particular during the famine that struck Russia in 1841. They were therefore no longer able to help him, and his salary as an adjunct professor was modest. The above mentioned biography states that for that reason, Chebyshev became very thrifty and remained mso throughout his life, but the only thing for which he never spared money was the materials he needed to construct his mechanical models. He was able to spend hundreds of thousands of roubles for his machines. From 1847 onward, Chebyshev was a university professor in Saint Petersburg, and thus although he was not wealthy, his income allowed him to live decently, and he was in no need of extra income. He worked for the army (with or without compensation) on his own accord.

\begin{itemize}

\item[$\bullet$]We conclude that Chebyshev was not in need of extra money.

\end{itemize}

{\bf The Crimean War.} There is no evidence of Chebyshev having anything to do with uniforms at the time of the Crimean War. No biography of Chebyshev mentions his working in the garment industry or in any other capacity related to military uniforms. His wartime work was connected to the artillery, probably in ballistics. It was common in the eighteenth and nineteenth centuries for mathematicians to perform work for the army and to teach in military schools. In France, Poisson taught at Saint Cyr military academy, a school for army officers, and like Monge, Carnot, Liouville, Hermite, Jordan, and many others, he also taught at the \'Ecole Polytechnique, which was (and still is in some sense) a military school.

\begin{itemize}

\item[$\bullet$]The Crimean War lasted from 1853 to 1856, ending more than twenty years before the 1878 article by Chebyshev on garment cutting.

\end{itemize}

{\bf The Book.} There is no book published by Chebyshev on the subject. Grant mentions the existence of a book inRussian that was never translated into a Western European language.

\begin{itemize}

\item[$\bullet$]It appears that such a book simply does not exist. It never existed.

\end{itemize}

The best summary of the situation can be found in the Dictionary of Scientific Biography \cite{Youschkevitch1970}:

\begin{itemize}

\item[]{\it From 1856 Chebyshev was a member of the Artillery Committee, which was charged with the task of introducing artillery innovations into the Russian army. In close cooperation with the most eminent Russian specialists in ballistics, such as N. V. Maievski, Chebyshev elaborated mathematical devices for solving artillery problems. He suggested (1867) a formula for the range of spherical missiles with initial velocities within a certain limit; this formula was in close agreement with experiments. Some of Chebyshev's works on the theory of interpolation were the result of the calculation of a table of fire effect based on experimental data. Generally, he contributed significantly to ballistics.}

\item[]{\it His paper of 1878, {\sc Sur la coupe des v\^etements}, provided the basis for a new branch of the theory of surfaces. Chebyshev investigated a problem of binding surface with cloth that is formed in the initial flat position with two systems of nonextensible rectilinear threads normal to one another. When the surface is bound with cloth the {\sc Chebyshev net}, whose two systems of lines form curvilinear quadrangles with equal opposite sides, appears. Wrapping a surface in cloth is a more general geometrical transformation than is deformation of a surface, which preserves the lengths of all the curved lines; distances between the points of the wrapped cloth that are situated on different threads are, generally speaking, changed in wrapping. In recent decades Chebyshev's theory of nets has become the object of numerous studies.}

\end{itemize}

\section*{Conclusions}

In his article of 2011, Papadopoulos expressed doubts concerning the authenticity of what Grant asserted in his contribution to MacTutor. This conclusion is confirmed by our previous analysis. The origin of the interest of Chebyshev on garment cutting is the one expressed by himself in his article of 1878. All the evidence indicates that the contribution by Grant on MacTutor is a fabrication. The most important lesson we may learn from this episode is the peril of citing sources without checking their veracity, for otherwise, we may be guilty of perpetuating inaccuracies and falsehoods.

\section*{Acknowledgments}

The content of this article grew from communications with Professors Athanese Papadopoulos and Galina Sinkevich. The author thanks Professor Denis Monod-Broca, secretary general of the {\it Association fran\c caise pour l'avancement des sciences} (AFAS), for giving permission to reproduce the pages from {\it C. R. Assoc. Fran\c caise l'Avanc. Sci}. This article was much enriched by commentaries and information provided by an anonymous referee. Stylistic remarks by managing editor David Kramer are greatly appreciated.


\end{document}